\documentclass[reqno]{amsart}

\usepackage{amsmath,mathbbol,amssymb,mathbbol}
\usepackage{epic}
\usepackage{pstricks}
\usepackage{amsthm}

\newtheorem{theorem}{Theorem}

\newtheorem{lemma}[theorem]{Lemma}

\theoremstyle{definition}

\newtheorem{remark}[theorem]{Remark}

\newtheorem*{thA}{{Theorem A}}
\newtheorem*{thB}{{Theorem B}}
\newtheorem*{thC}{{Theorem C}}

\def\R{\mathbb{R}}

\def\N{\mathbb{N}}
\def\C{\mathbb{C}}
\def\ZZ{\mathbb{Z}}

\def\E{\mathbb{E}}

\def\var{\varepsilon}

\def\l{\lambda}

\def \supp{\operatorname{supp}}

\def\cC{{\mathcal C}}
\def\cD{{\mathcal D}}
\def\cE{{\mathcal E}}

\def\cL{{\mathcal L}}

\def\cR{{\mathcal R}}

\def\cV{{\mathcal V}}
\def\cW{{\mathcal W}}
\def\cI{{\mathcal I}}
\def\cJ{{\mathcal J}}

\def\var{\textrm{Var}}

\title[Fourier multipliers on weighted $L^p$ spaces]
{Fourier multipliers on weighted $L^p$ spaces}

\author{Sebastian Kr\'ol}
\address{S. Kr\'ol, Faculty of Mathematics and Computer Science, Nicolaus Copernicus University, ul. Chopina 12/18, 87-100 Toru\'n, Poland, and Institut f\"ur Analysis, Fachrichtung Mathematik, TU Dresden, 01062 Dresden, Germany}
\email{sebastian.krol@mat.umk.pl}

\begin{document}

\thanks{The author was supported by the Alexander von Humboldt Foundation}



\begin{abstract}
The paper provides a complement to the classical results
on Fourier multipliers on $L^p$ spaces.
In particular, we prove that if $q\in (1,2)$ and a function
$m:\R \rightarrow \C$ is of bounded $q$-variation uniformly on the dyadic intervals in $\R$, i.e. $m\in V_q(\cD)$, then $m$ is a Fourier multiplier on $L^p(\R, wdx)$ for every $p\geq q$ and every weight $w$  satisfying Muckenhoupt's $A_{p/q}$-condition. We also obtain a higher dimensional counterpart of this result as well as of a result by E. Berkson and T.A.~ Gillespie including the case of the $V_q(\cD)$ spaces with $q>2$.
New weighted estimates for modified Littlewood-Paley functions are also provided.\\[-8ex]
\end{abstract}

\renewcommand{\subjclassname}{\textup{2010} Mathematics Subject Classification}

\maketitle

\section{Introduction and Statement of Results}

For an interval $[a,b]$ in $\R$ and a number $q\in [1,\infty)$ denote by
$V_q([a,b])$ the space of all functions $m:[a,b] \rightarrow \C$ of bounded $q$-variation
 over $[a,b]$, i.e.,
$$\|m\|_{V_q([a,b])}:= \sup_{x\in [a,b]}|m(x)| + \|m\|_{\var_q([a,b])} <\infty,
$$
where $\|m\|_{\var_q([a,b])}:=\sup\{ (\sum_{i=0}^{n-1} |m(t_{i+1}) - m(t_i)|^q )^{1/q}\}$ and the supremum is taken over all finite sequences $a=:t_0< t_1<...< t_n := b$  $(n\in \N)$.
We write $\cD$ for the dyadic decomposition of $\R$, i.e.,
$\cD:=\left\{ \pm(2^k,2^{k+1}]: k\in \ZZ \right\}$, and set
 $$V_q(\cD):=\left\{m:\R\rightarrow \C :\quad \sup_{I\in \cD} \|m_{|I}\|_{V_q(I)} <\infty \right\} \quad \quad (q\in [1,\infty)).$$

Moreover, let $A_p(\R)$ $(p\in [1,\infty))$ be the class of weights on $\R$ which satisfy the Muckenhoupt $A_p$ condition. Denote by $[w]_{A_p}$ the $A_p$-constant of $w\in A_p(\R)$.
If $w \in A_\infty(\R):=\cup_{p\geq 1} A_p(\R)$ we write $M_p(\R,w)$ for the class of all  multipliers
 on $L^p(\R, w)$ $(p> 1)$, i.e.,
$$M_p(\R, w) := \left\{ m\in L^\infty(\R): T_m \textrm{ extends
 to a bounded operator on } L^p(\R, w) \right\}.$$
Here $T_m$ stands for the Fourier multiplier
with the symbol $m$, i.e., $(T_mf) \widehat{} = m \widehat{f}$ $(f\in S(\R)).$
Note that $M_p(\R,w)$ becomes a Banach space under the norm
$\|m\|_{M_p(\R,w)} := \|T_m\|_{\cL(L^p(\R, w))}$ $(m\in M_p(\R, w) )$.

The main result of the paper is the following complement to results due to D.~Kurtz \cite{Kurtz80},
R. Coifman, J.-L. Rubio de Francia, S. Semmes \cite{CRdeFS88}, and  E. Berkson,
 T. Gillespie \cite{BeGi98}.

\begin{thA}\label{main}\emph{
$(i)$ Let $q\in (1,2]$. Then, $V_q(\cD)\subset M_p(\R,w)$ for every
$p\geq q$ and every Muckenhoupt weight $w\in A_{p/q}(\R)$. \\
$(ii)$ Let $q>2$. Then, $V_q(\cD)\subset M_p(\R, w)$ for every
$2\leq p <(\frac{1}{2} - \frac{1}{q})^{-1}$ and every Muckenhoupt weight $w\in A_{p/2}$ with
$s_w> (1- p(\frac{1}{2} - \frac{1}{q}))^{-1}$.}
\end{thA}
Here, for every $w\in A_\infty(\R)$, we set $s_w:=\sup\{ s\geq 1: w \in RH_s(\R) \}$ and we write $w\in RH_s(\R)$ if $$\sup_{a<b} \left( \frac{1}{b-a}\int_a^b w(x)^s dx \right)^{1/s}
\left(\frac{1}{b-a} \int_a^b w(x) dx\right)^{-1} <\infty.$$
Recall that, by the reverse H\"older inequality,  $s_w \in (1, \infty]$ for every Muckenhoupt weight $w\in A_\infty(\R)$.

For the convenience of the reader we repeat the relevant material from the literature, which we also use in the sequel.

Recall first that in \cite{Kurtz80} D. Kurtz proved the following weighted variant of the
classical Marcinkiewicz multiplier theorem.

\begin{theorem}[{\cite[Theorem 2]{Kurtz80}}]\label{Kurtz} $V_1(\cD) \subset M_p(\R,w)$ for every
$p\in (1,\infty)$ and every Muckenhoupt weight $w\in A_p(\R)$.
\end{theorem}

As in the unweighted case, Theorem \ref{Kurtz} is equivalent to a weighted variant of
the Littlewood-Paley decomposition theorem, which asserts that for
the square function $S^\cD$ corresponding to the dyadic decomposition $\cD$ of $\R$,
$\|S^\cD f\|_{p,w}\eqsim \|f\|_{p,w}$ $(f\in L^p(\R,w))$
 for every $p\in (1,\infty)$ and $w\in A_p(\R)$; see \cite[Theorem 1]{Kurtz80}, and also \cite[Theorem 3.3]{Kurtz80}.
Here and subsequently, if $\cI$ is a family of disjoint intervals in $\R$, we write $S^\cI$ for the Littlewood-Paley square function corresponding to $\cI$, i.e.,
$S^\cI f : = \left(\sum_{I\in \cI} |S_If|^2 \right)^{1/2}$
 $(f\in L^2(\R)).$

Recall also that in \cite{RdeF85} J.-L. Rubio de Francia proved the following
extension of the classical Littlewood-Paley decomposition  theorem.

\begin{theorem}[{\cite[Theorem 6.1]{RdeF85}}]\label{RdeFtheorem}
Let $2< p <\infty$ and $w\in A_{p/2}(\R)$. Then for an arbitrary family $\cI$ of disjoint  intervals in $\R$ the square function $S^\cI$ is bounded on $L^p(\R, wdx)$.
\end{theorem}

Applying Rubio de Francia's inequalities, i.e. Theorem \ref{RdeFtheorem},  R.
Coifman, J.-L. Rubio de Francia, and S. Semmes \cite{CRdeFS88} proved the following
extension and improvement of the classical Marcinkiewicz multiplier theorem.
(See Section 2 for the definition of $R_2(\cD)$.)

\begin{theorem}[{\cite[Th\'eor\`eme 1 and Lemme 5]{CRdeFS88}}] \label{CRdeFS}
Let $2\leq q<\infty$. Then, $V_q(\cD)\subset M_p(\R)$ for every $p\in (1,\infty)$ such that
$|\frac{1}{p} - \frac{1}{2}|<\frac{1}{q}$.

Furthermore, $R_2(\cD) \subset M_2(\R, w)$ for every $w\in A_1(\R)$.
\end{theorem}

 Subsequently, a weighted variant of Theorem \ref{CRdeFS} was given by E. Berkson and T.
Gillespie in \cite{BeGi98}. According to our notation their result can be formulated as follows.

\begin{theorem}[{\cite[Theorem 1.2]{BeGi98}}] \label{BeGi}
Suppose that $2\leq p<\infty$ and $w\in A_{p/2}(\R)$. Then,
 there is a real number $s>2$, depending only on $p$ and $[w]_{A_{p/2}}$, such that $\frac{1}{s}>\frac{1}{2} - \frac{1}{p}$ and $V_q(\cD) \subset M_p(\R, w)$ for all
$1\leq q <s$.
\end{theorem}

Note that the part $(i)$ of Theorem A fills a gap which occurs in Theorem \ref{Kurtz} and the weighted part of Theorem \ref{CRdeFS}. The part $(ii)$ identifies the constant $s$ in  Berkson-Gillespie's result, i.e., Theorem \ref{BeGi}, as $(\frac{1}{2} - \frac{1}{s'_w p})^{-1}$, where $s'_w:=\frac{s_w}{s_w-1}$, and in general, this constant is  best possible.

Except for some details, the proofs given below reproduce well-known arguments from the Littlewood-Paley theory; in particular, ideas which have been presented in \cite{Kurtz80}, \cite{CRdeFS88}, \cite{RdeF85}, and \cite{Xu96}.
A new point of our approach is the following result on weighted estimates for modified Littlewood-Paley functions $S^\cI_q(\cdot) :=(\sum_{I\in\cI} |S_I (\cdot)|^{q'})^{1/q'}$ $(q\in   (1,2])$, which may be of independent interest. 

\begin{thB}\label{LPineq}{\it{
 $(i)$ Let $q\in (1,2)$, $p>q$, and $w\in A_{p/q}(\R)$.
 Then, there exists a constant $C>0$ such that for any 
 family $\cI$ of disjoint intervals in $\R$ 
 $$\|S^\cI_q f\|_{p,w} \leq C\|f\|_{p,w}  \quad\quad (f\in L^p(\R, wdx)).$$
Moreover, for every $q\in (1,2)$, $p>q$ and $\cV \subset A_{p/q}(\R)$ with $\sup_{w\in \cV}[w]_{A_{p/q}}<\infty$  
$$\sup\left\{ \|S^\cI_qf\|_{p,w}: w\in \cV, \; \cI \textrm{ a family of disjoint intervals in }\R, \; \|f\|_{p,w}=1 \right\}
<\infty.$$

$(ii)$ For any family $\cI$ of disjoint intervals in $\R$  and every Muckenhoupt weight $w \in A_{1}(\R)$, the operator $S^\cI_2$ maps $L^2(\R,wdx)$ into weak-$L^2(\R,wdx)$, and 
$$\sup\left\{ \|S^\cI_2f\|_{L_w^{2,\infty}}: w\in \cV, \; \cI \textrm{ a family of disjoint intervals in }\R, \|f\|_{L_w^{2}}=1\right\}<\infty$$
for every  $\cV \subset A_{1}(\R)$ with $\sup_{w\in \cV}[w]_{A_1}<\infty$. 

Moreover, if $q\in (1,2)$, then for any well-distributed family $\cI$ of disjoint intervals in $\R$ and every Muckenhoupt weight $w\in A_1(\R)$, the operator $S_q^\cI$ maps $L^q(\R,wdx)$ into weak-$L^q(\R,wdx)$. }}
\end{thB}
Recall that a family $\cI$ of disjoint intervals in $\R$ is {\it well-distributed} if there exists $\lambda>1$ such that $\sup_{x\in \R}\sum_{I\in\cI} \chi_{\lambda I}(x) < \infty$, where $\lambda I$ denotes the interval with the same center as $I$ and length $\lambda$ times that of $I$.

Note that the validity of the $A_1$-weighted $L^2$-estimates for square function $S^\cI=S^\cI_2$ corresponding to an arbitrary family $\cI$ of disjoint intervals in $\R$, i.e., 
$$\|S^\cI_2 f\|_{2,w} \leq C_{w}\|f\|_{2,w}  \quad\quad (f\in L^2(\R, wdx), \; w\in A_1(\R)),$$ 
is conjectured by J.-L. Rubio de Francia in \cite[Section 6, p.10]{RdeF85}; see also 
\cite[Section 8.2, p. 187]{Duo01}. Theorem B$(ii)$, in particular,
provides the validity of the weak variant of Rubio de Francia's conjecture.
Notice that in contrast to the square function operators $S_2^\cI$, in general, operators $S^\cI_q$ $(q\in[1,2))$ are not bounded on (unweighted) $L^{q}(\R)$; see \cite{CoTa01}.
Moreover, in \cite{Quek03} T.S. Quek proved that if $\cI$ is a well-distributed family of disjoint intervals in $\R$, then the operator $S^\cI_q$ maps $L^q(\R)$ into
$L^{q,q'}(\R)$ for every $q\in (1,2)$. Note that this result is in a sense sharp, i.e., 
$L^{q, q'}(\R)$ cannot be replaced by $L^{q,s}(\R)$ for any $s<q'$; see \cite[Remark 3.2]{Quek03}.  Therefore, Theorem B provides also a weighted variant of this line of researches. 
Cf. also relevant results given by S.V. Kisliakov in \cite{Kis08}.

Furthermore, as a consequence of our approach we also get a higher dimensional analogue of Theorem A, see Theorem C in Section 4, which extends earlier results by Q. Xu \cite{Xu96}; see also M. Lacey \cite[Chapter 4]{Lac07}.
Since the formulation of Theorem~C is more involved and its proof is essentially the iteration of one-dimensional arguments we refer the reader to Section 4 for more information.

The part $(ii)$ of Theorem A is a quantitative improvement of 
\cite[Theorem 1.2]{BeGi98} due to E. Berkson and T. Gillespie. Furthermore, we present an alternative approach based on a version of the Rubio de Francia extrapolation theorem that holds for limited ranges of $p$ which was recently given in \cite{AuMa07}.

The organisation of the paper is well-reflected by the titles of the following sections. However, we conclude with an additional comment.
The proof of Theorem~A is based on weighted estimates from the part $(i)$ of Theorem B.
To keep the pattern of the proof of the main result of the paper, Theorem A, more transparent, we postpone the proof of Theorem B$(ii)$ to Section 3.

\section{Proofs of Theorems B$(i)$ and A}

We first introduce auxiliary spaces which are useful in the proof of Theorem~A.
Let $q\in[1,\infty)$. If $I$ is an interval in $\R$ we denote by $\cE(I)$
the family of all step functions from $I$ into $\C$.
If $m:=\sum_{J\in \cI} a_J \chi_J$, where $\cI$ is a decomposition of $I$ into subintervals and $(a_J)\subset \C$, write $[m]_{q}:=(\sum_{J\in \cI} |a_J|^q )^{1/q}$. Set
$\cR_q(I):=\left\{ m \in \cE(I) : [m]_q\leq 1 \right\}$
and
$$\cR_q(\cD):= \left\{ m:\R \rightarrow \C:\,\, m_{|I} \in \cR_q(I)\,\,\textrm{for every } I\in \cD \right\}.$$
Moreover, let
$$R_q(I):=\left\{ \sum_{j} \l_j m_j : \, m_j \in \cR_q(I), \,\, \sum_j |\l_j|<\infty    \right\}$$ and
 $$\|m\|_{R_q(I)} :=\inf\left\{\sum_j |\l_j| :\, m=\sum_j \l_j m_j, \,\, m_j\in \cR_q(I) \right\}\quad \left(m\in R_q(I) \right).$$
Note that $\left(R_q(I), \|\cdot\|_{R_q(I)}\right)$ is a Banach space.
Set
$$R_q(\cD):=\left\{m:\R \rightarrow \C:\, \sup_{I\in \cD} \|m_{|I}\|_{R_q(I)} <\infty \right\} \quad \quad (q\in [1,\infty)).$$

In the sequel, if $\cI$ is a family of disjoint intervals in $\R$, we write $S^\cI_1 f:=\sup_{I\in \cI}|S_I f|$ ($f\in L^1(\R)$) and
 $S^\cI_r f :=(\sum_{I\in\cI} |S_I (f)|^{r'})^{1/r'}$ $(r\in (1,2], f\in L^r(\R))$.

We next collect main ingredients of the proof of
Theorem B$(i)$, which provides crucial vector-valued estimates for weighted multipliers 
in the proof of Theorem A; see e.g. \eqref{vectorestim}.

Lemma \ref{CarHunt} is a special version of the result on weighted inequalities for Carleson's operator given by J.-L. Rubio de Francia, F. J. Ruiz and J. L. Torrea in \cite{RuRuTo86}; see also \cite[Remarks 2.2, Part III]{RuRuTo86}.

\begin{lemma}[{\cite[Theorem 2.1, Part III]{RuRuTo86}}]\label{CarHunt} 
Let $s\in (1,\infty)$ and $w\in A_s(\R)$.
Then, there exists a constant $C>0$ such that for any family $\cI$ of disjoint intervals in 
$\R$
$$\|S^\cI_1 f\|_{s,w} \leq C\|f\|_{s,w}  \quad (f\in L^s(\R, wdx)).$$

Moreover, for every $s>1$ and every set $\cV \subset A_s(\R)$ with $\sup_{w\in \cV}[w]_{A_s}<\infty$ 
$$\sup\left\{ \|S^\cI_1\|_{s,w}: w\in\cV,\;  \cI \textrm{ a family of disjoint intervals in }\R \right\} <\infty.$$
\end{lemma}

\begin{remark}\label{upperbound}
The second statement of Lemma \ref{CarHunt} can be obtained from a detailed analysis of the constants involved in the results which are used in the proof of \cite[Theorem 2.1$(a)\Rightarrow (b)$, Part III]{RuRuTo86}, i.e., the weighted version of the Fefferman-Stein inequality and the reverse H\"older inequality.

Recall the weighted version of the Fefferman-Stein inequality, which in particular says that for every $p\in (1,\infty)$ and every Muckenhoupt weight $w\in A_p(\R)$ there exists a constant $C_{p,w} >0$, which depends only on $p$ and $[w]_{A_p}$, such that
\begin{equation}\label{FS}
\int_{\R} M f(t)^p \, w(t) \; dt\leq C_{p,w} \int_{\R}
M^\# f(t)^p \, w(t) \; dt  \quad (f\in L^p(\R)\cap L^p(\R,w) ),
\end{equation}
where $M$ and $M^\sharp$ denote the Hardy-Littlewood maximal operator and the Fefferman-Stein sharp maximal operator, respectively; see \cite[Theorem, p.41]{Jo83}, or \cite[Theorem 2.20, Chapter IV]{GCRu85}. 
We emphasize here that the constant $C_{p,w}$ on the right-hand side of this inequality is not given explicitly in the literature, but it can be obtained from a detailed analysis of the constants involved in the results which are used in the proof of \eqref{FS}, 
$\sup_{w\in\cV} C_{p,w} <\infty$ for every subset $\cV\subset A_p(\R)$ with 
$\sup_{w\in \cV}[w]_{A_p}<\infty$.

Furthermore, it should be noted that if $\cV \subset A_p(\R)$  with 
$\sup_{w\in \cV}[w]_{A_p}<\infty$, then there exists $\epsilon>0$ such that $\cV \subset A_{p-\epsilon}(\R)$ and $\sup_{w\in \cV}[w]_{A_{p-\epsilon}} <\infty$.
 It can be directly obtained from a detailed analysis of the constants involved in main ingredients of the proof of the reverse H\"older inequality. Cf., e.g., \cite[Lemma 2.3]{LeOm12}. 

We refer the reader to \cite[Chapter IV]{GCRu85} and \cite[Chapter 7]{Duo01} for recent expositions of the results involved in the proof of the reverse H\"older inequality and the Fefferman-Stein inequality, which originally come from \cite{CoFe74}, and \cite{Mu72}, \cite{Mu74}. 

\end{remark}

The next lemma is a special variant of Rubio de  Francia's extrapolation theorem; see \cite[Theorem 3]{RdeF85}. For the convenience of the reader we rephrase \cite[Theorem 3]{RdeF85} here in the context of Muckenhoupt weights merely.

\begin{lemma}[{\cite[Theorem 3]{RdeF84}}]\label{extrapol}
Let $\lambda$ and $r$ be fixed with
$1\leq \lambda \leq r <\infty$, and let $\mathcal{S}$ be a family of
 sublinear operators which is uniformly bounded in $L^r(\R,wdx)$ for each
$w\in A_{r/\lambda}(\R)$, i.e.,
$$\int|Sf|^r wdx \leq C_{r,w} \int |f|^rwdx \quad (S\in \mathcal{S}, \,
w\in A_{r/\lambda}(\R)).$$
If $\lambda <p, \alpha < \infty$ and $w \in A_{p/\lambda}(\R)$, then $\mathcal{S}$
is uniformly bounded in $L^p(\R, wdx)$ and
even more:
\begin{equation*}\label{eq1extr}
\int(\sum_{j}|S_jf_j|^\alpha)^{p/\alpha} wdx \leq C_{p,\alpha,w}
 \int (\sum_{j} |f_j|^\alpha)^{p/\alpha} wdx\quad \left( f_j\in L^p(\R, wdx),\,\, S_j\in \mathcal{S} \right).
\end{equation*}
\end{lemma}

Combining Lemma \ref{CarHunt} with Theorem \ref{RdeFtheorem} we get the intermediate weighted estimates for operators $S_q^\cI$ $(q\in(1,2))$ stated in Theorem B$(i)$.

For the background on the interpolation theory we refer the
reader to \cite{BeLo77}; in particular, see \cite[Chapter 4 and Section 5.5]{BeLo77}.

\begin{proof}[{\bf{Proof of Theorem B$(i)$}}]
Fix $q\in(1,2)$ and $w\in A_{2/q}(\R)$. By the reverse H\"older inequality, 
$w\in A_{2/r}(\R)$ for some $r\in(q,2)$. 
Note that there exist $p\in(2,q')$ and $s>1$ such that $\frac{p}{q'} \frac{1}{p} + 
(1 - \frac{p}{q'})\frac{1}{s}  = \frac{1}{r}$.
Therefore, combining Theorem \ref{RdeFtheorem} with Lemma \ref{CarHunt}, by complex interpolation, the operator $S^{\cI}_{({2q'}/{p})'}$ is bounded on $L^r(\R, v)$ for every $v\in A_1(\R)$. Since $p>2$, the same conclusion holds for $S^\cI_{q}$.

By Rubio de Francia's extrapolation theorem, Lemma \ref{extrapol}, we get that $S^\cI_{q}$ is bounded on $L^2(\R, v)$ for every $v\in A_{2/r}(\R)$. According to our choice of $r$, we get the boundedness of $S^\cI_{q}$ on $L^2(\R, w)$.

Since the weight $w$ was taken arbitrarily, we can again apply Rubio de Francia's extrapolation theorem, Lemma \ref{extrapol}, to complete the proof of the first statement.

The second statement follows easily from a detailed analysis of the first one. For a discussion on the character of the dependence of constants in Rubio de Francia's iteration algorithm, we refer the reader to \cite{DGPP05}, or \cite[Section 3.4]{CrMaPe11}. See also the comment on 
the reverse H\"older inequality in Remark \ref{upperbound}.
\end{proof}

Note that
$R_q(I) \varsubsetneq V_q(I)$  for every
interval $I$ in $\R$ and $q\in [1,\infty)$.
However, the following reverse inclusions hold for these classes.

\begin{lemma}[{\cite[Lemme 2]{CRdeFS88}}]\label{C-RdeF-S}
 Let $1\leq q< p<\infty$. For every interval $I$ in $\R$,
$V_q(I) \subset R_p(I)$ with the inclusion norm bounded by a constant independent of $I$.
\end{lemma}

The patterns of the proofs of the parts $(i)$ and $(ii)$ of Theorem A are essentially the same.
Therefore, we sketch the proof of the part $(ii)$ below.

\begin{proof}[{\bf{Proof of Theorem A}}]
$(i)$ We only give the proof for the more involved case $q\in (1,2)$; the case $q=2$ follows simply from Theorem \ref{CRdeFS} and interpolation arguments presented below; see also Remark \ref{rems3} below. 

Fix $q\in (1,2)$. We first show that for every subset $\cV\subset A_{{2}/{q}}(\R)$ such that $\sup_{w\in \cV} [w]_{A_{{2}/{q}}} <\infty$ we have
$$\sup\left\{ \|T_{m\chi_I}\|_{2,w}: m\in R_q(\cD), \|m\|_{R_q(\cD)}\leq 1,  w\in \cV, I\in \cD  \right\} <\infty.$$

Fix $\cV\subset A_{{2}/{q}}(\R)$ with $\sup_{w\in \cV} [w]_{A_{{2}/{q}}}<\infty$. Note that, by the definition of the $R_q$-classes, it is sufficient to prove 
the claim with $R_q(\cD)$ replaced by $\cR_q(\cD)$ .
Fix $m\in \cR_q(\cD)$ and set
 $ m\chi_I =:\sum_{J\in \cI_{I}} a_{I,J} \chi_J$ for every $I\in \cD$,
where $\cI_I=\cI_{I,m}$ is a decomposition of $I$ and
$(a_{I,J})_{J\in \cI_{I}}\subset \C$ is a sequence with $\sum_{J\in \cI_{I}}|a_{I,J}|^q\leq 1$.
Note that $T_{m\chi_I}f = \sum_{J} a_{I,J} S_Jf$ and $\|T_{m\chi_I}f\|_{2,w} \leq \|S_q^{\cI_I}f\|_{2,w}$ for every $I\in \cD$, $w\in \cV$ and $f\in L^{2}(\R,w)$. Therefore, by Lemma \ref{CarHunt}, our claim holds. 

By interpolation argument, we next sharpen this claim and prove that for every subset
$\cV\subset A_{2/q}(\R)$ with $\sup_{w\in\cV}[w]_{A_{2/q}}<\infty$ there exists 
$\alpha=\alpha(q, \cV)>1$ such that 
\begin{equation}\label{step1}
\sup\left\{ \|T_{m\chi_I}\|_{2,w}: m\in R_{\alpha q}(\cD), \|m\|_{R_{\alpha q}(\cD)}\leq 1,  w\in \cV, I\in \cD  \right\} <\infty.
\end{equation}

Note that, by the reverse H\"older inequality, see also Remark \ref{upperbound}, there exists $\alpha>1$ such that  $w^\alpha \in A_{2/q}(\R)$ $(w\in \cV)$ and $\sup_{w\in\cV}[w^\alpha]_{A_{2/q}}<\infty$.
From what has already been proved and Plancherel's theorem, for every $I\in \cD$ and $w\in \cV$ the bilinear operators
$$  R_q(I)\times L^{2}(\R, w^\alpha dx) \ni (m, f) \mapsto T_mf\in
 L^{2}(\R, w^\alpha dx) $$
$$  L^\infty(\R)\times L^{2}(\R)\ni (m, f) \mapsto T_mf\in
 L^{2}(\R) $$
are well-defined and bounded uniformly with respect to $w\in \cV$ and $I\in \cD$. Therefore, by complex interpolation,  $\left(R_q(I), L^\infty(\R)\right)_{[\frac{1}{\alpha}]} \subset M_{2}(\R, w)$.
However, it is easy to check that  $R_{\alpha q}(I)\subset \left(R_q(I), L^\infty(\R)\right)_{[\frac{1}{\alpha}]}$ with the inclusion norm bounded by a constant independent of $I\in \cD$. We thus get \eqref{step1}.

In consequence, by Lemma \ref{C-RdeF-S}, it follows that
\begin{equation}\label{vectorestim}
 \sup\left\{ \|T_{m\chi_I}\|_{2,w}: m\in V_{q}(\cD),\; \|m\|_{V_q(\cD)}\leq 1,\;  w\in \cV, I\in \cD  \right\} <\infty
\end{equation}
for every subset $\cV\subset A_{2/q}(\R)$ with $\sup_{w\in\cV}[w]_{A_{2/q}}<\infty$.

Hence, we can apply a truncation argument based on Kurtz' weighted variant of Littlewood-Paley's inequality. Namely, fix $w\in A_{2/q}(\R)$, $m\in V_q(\cD)$ with 
$\|m\|_{V_q(\cD)}\leq 1$, and 
 $f\in L^2(\R)\cap L^2(\R,w)$,  $g\in L^2(\R,w)\cap L^2(\R, w^{-1})$. 
Note that $gw\in L^2(\R)$ and $A_{2/q}(\R)\subset A_{2}(\R)$.
Therefore,  combining the Cauchy-Schwarz inequality and Kurtz' result, \cite[Theorem 1]{Kurtz80}, we get:
\begin{eqnarray*}
|(T_m f,g)_{L^2(\R,w)}| & = & \left|\sum_{I\in\cD} \int_{\R} S_I(T_mf) S_I(gw) dx\right| \\
& \leq  & C \left\| (\sum_{I\in\cD} |T_{m\chi_I}S_If|^2)^{1/2} \right\|_{2,w}
\left\| (\sum_{I\in\cD} |S_I (gw)|^2)^{1/2} \right\|_{2,w^{-1}}\\
& \leq & C \|f\|_{2,w} \|g\|_{2,w},
\end{eqnarray*}
where $C$ is an absolute constant independent of $m$, $f$ and $g$.
 Now the converse of H\"older inequality and a density argument show that
$m\in M_2(\R, w)$. 

Consequently,  $V_q(\cD)\subset M_2(\R, w)$, and Rubio de Francia's extrapolation theorem, 
Lemma \ref{extrapol}, yields $V_q(\cD) \subset M_p(\R, w)$ for every $p>q$ and every Muckenhoupt weight $w\in A_{p/q}(\R)$.

It remains to prove that $V_q(\cD) \subset M_q(\R, w)$ for every $w\in A_1(\R)$.
Fix $m\in V_q(\cD)$ and $w\in A_1(\R)$. Then, by Theorem \ref{CRdeFS} (see also
Remark \ref{rems3}), $T_m$ is bounded on $L^r(\R)$ for every $r\in (1,\infty)$. From what has already been proved, $T_m$ is bounded on $L^r(\R, w)$ for every $r>q$. Therefore, the boundedness of $T_m$ on $L^q(\R, w)$
follows by the reverse H\"older inequality for $w$ and a similar interpolation argument as before. This completes the proof of the part $(i)$.\\

$(ii)$ Fix $q>2$ and $s>\frac{q}{2}$. Let $\cV_{s}:=\{w\in A_1(\R): w\in RH_s(\R)\}$.
Note that there exists $r=r_s>q$ such that
$\frac{1}{s}\frac{1}{2}+\frac{1}{s'} \frac{1}{r} <\frac{1}{q}.$

Fix $w \in \cV_{s}$. By Theorem \ref{CRdeFS}, the bilinear operators
$$R_r(\cD) \times L^2(\R)\ni (m,f) \mapsto T_m f \in L^2(\R)$$
$$R_2(\cD) \times L^2(\R, w^{s}) \ni (m,f) \mapsto T_mf \in L^2(\R, w^{s})$$
are well-defined and bounded.
By interpolation, it follows that
$$M_2(\R, w) \supset  \left(R_2(\cD), R_r(\cD)\right)_{[\frac{1}{s}]}
\supset R_{\alpha q}(I) $$ uniformly with respect to $I\in \cD$, where
$\alpha= \alpha_s:= (\frac{1}{2s} + \frac{1}{s' r})^{-1}/q>1$.

As in the corresponding part of the proof of $(i)$, by
 truncation and duality arguments,
 we get $R_{\alpha q}(\cD)\subset M_2(\R, w).$

Consequently, since $\alpha_s>1$ for every $s>\frac{q}{2}$, by Lemma \ref{C-RdeF-S},
\begin{equation}
 V_q(\cD)\subset M_2(\R, w)\quad \quad \textrm{for every }
w \in \bigcup_{s>\frac{q}{2}} \cV_{s}(\R).
\end{equation}
Note that this is precisely the assertion of $(ii)$ for $p=2$.

We can now proceed by extrapolation. Since for every $s>\frac{q}{2}$ we can rephrase
$\cV_{s}$ as $A_{\frac{2}{2}}(\R)\cap RH_{(\frac{2s'}{2})'}(\R)$,
by \cite[Theorem 3.31]{CrMaPe11}, we get
\begin{equation}\label{pextra}
 V_q(\cD)\subset M_p(\R, w)\quad \textrm{for every }
 s>\frac{q}{2}, \,\,2<p < 2s',  \textrm{ and }
 w \in A_{\frac{p}{2}}(\R)\cap RH_{(\frac{2s'}{p})'}(\R)\,.
\end{equation}
Finally, it is easy to see that for every
$2\leq p < \frac{1}{2} - \frac{1}{q}= 2(\frac{q}{2})'$ and $w \in A_{p/2}(\R)$ with
$s_w >(1-p(\frac{1}{2} - \frac{1}{q}))^{-1} = (\frac{2}{p}(\frac{q}{2})')'$ there exists
$s=s_{p,w} >\frac{q}{2}$ such that $p<2s'$ and $w \in RH_{(\frac{2s'}{p})'}$.
Therefore, \eqref{pextra} completes the proof of $(ii)$.
\end{proof}

\begin{remark}\label{rems3}
 In the proof of Theorem A  we use Theorem \ref{CRdeFS} due to Coifman, Rubio de Francia, Semmes. Note that the patterns of all proofs are essentially the same.

Indeed, we can rephrase the proof of \cite[Th\'eor\`eme 1]{CRdeFS88} as follows.
First recall that $M_p(\R)=M_{p'}(\R)$ for every $p\in (1,\infty)$.
Let $r\geq 2$. By the Littlewood-Paley decomposition theorem,  Rubio de Francia's
inequalities, and Plancherel's theorem,  the bilinear operators
$$ R_2(\cD) \times L^r(\R)\ni (m, f)\mapsto T_m f\in L^r(\R)$$
$$ L^\infty(\R)\times L^2(\R) \ni (m,f) \mapsto T_m f \in L^2(\R)$$
are well defined and bounded. Therefore, by interpolation,
$\left(R_2(\cD),L^\infty(\R)\right)_{[\theta(r)]} \subset M_p(\R)$, where
$\theta(r)\in (0, 1)$ and $p$ such that $\frac{1}{p}=\theta(r) \frac{1}{r} + (1-\theta(r)) \frac{1}{2}$.

Note that if $p\geq 2$ and $q$ satisfies $\frac{1}{q}>\frac{1}{2} -
 \frac{1}{p}$, then there exists $r>2$ such that
$R_{\alpha q}(I) \subset \left(R_2(\cD),L^\infty(\R)\right)_{[\theta(r)]}$ for an appropriate $\alpha>1$ and uniformly with respect to $I\in \cD$. Indeed, $\frac{1}{2}\theta(r) \searrow \frac{1}{2} - \frac{1}{p}$
 as $r\rightarrow \infty$. Therefore, Lemma \ref{C-RdeF-S} completes the
proof of Theorem \ref{CRdeFS}$(i)$.
\end{remark}

\section{Proof of Theorem B$(ii)$}

We obtain the proof of Theorem B$(ii)$ by means of a Banach function space analogue of Kurtz' weighted variant of Littlewood-Paley inequalities and the Fefferman-Stein inequality; see Lemma \ref{wellextr} below.

Note that without loss of generality in the proof of Theorem B$(ii)$ one can consider only
families consisting of bounded intervals in $\R$. 
For a bounded interval $I\in \cI$ we write  $\cW_I$ for Whitney's decomposition of  $I$ (see \cite[Section 2]{RdeF85} for the definition).
Note also that each decomposition $\cW_I$, $i\in \cI$, is of dyadic type. Furthermore, the family $\cW^\cI:=\bigcup_{I\in \cI} \cW_I$ is well-distributed, i.e.,
$$\sup_{x\in \R} \sum_{I\in \cW^\cI} \chi_{2I}(x)\leq 5.$$

We refer the reader primarily to \cite{BeSh88} for the background on function spaces.
In the sequel, let $\E$ denote a rearrangement invariant Banach function space over $(\R,dx)$.
Recall that, by Luxemburg's representation theorem \cite[Theorem 4.10, p.62]{BeSh88}, there exists a rearrangement invariant Banach function space $\overline{\E}$ over ($\R_+, dt$) such that for every
scalar, measurable function $f$ on $\R$,
$f\in \E$ if and only if $f^*\in \overline\E$, where $f^*$ stands for the decreasing rearrangement of $f$.
In this case $\|f\|_\E=\|f^*\|_{\overline\E}$ for every $f\in \E$.

Following \cite{LiTz79}, we define the \it {lower} \rm and \it {upper Boyd indices} \rm respectively by
$$
p_\E  := \lim_{t\to \infty} \frac{\log t}{\log h_\E(t)} \quad\quad \textrm{and}\quad\quad
q_\E  := \lim_{t\to 0+} \frac{\log t}{\log h_\E(t)}, 
$$
where $h_{\E}(t)=\|D_t\|_{\cL (\overline\E)}$ and
$D_t: \overline\E\rightarrow\overline\E$ $(t>0)$ is the \it {dilation operator} \rm defined by
$$
D_tf(s)=f(s/t), \qquad 0<t<\infty, \quad f\in \overline\E.
$$
One always has $1\le p_\E\le q_\E\le\infty$, see for example \cite[Proposition 5.13, p.149]{BeSh88}, where the Boyd indices are defined as the reciprocals with respect to our definitions.

Let $w$ be a weight in $A_\infty(\R)$.
Then we can associate with $\E$ and $w$ a rearrangement invariant Banach function space over $(\R, wdx)$ as follows
$$
\E_w =\{ f:\R \rightarrow \C \;\; {\rm{measurable }}:  f^*_w \in \overline{\E} \},
$$
and its norm is $\|f\|_{\E_w}=\| f^*_w\|_{\overline\E}$, where $f^*_w$ denotes the decreasing rearrangement of $f$ with respect to $wdx$.

For further purposes, recall also that examples of rearrangement Banach function spaces are the Lorentz spaces $L^{p,q}$ ($1\leq p$, $q\leq \infty$).
Note that $L^{p,\infty}_w=\textrm{weak}-L^p(\R, w)$ for every $p\in(1,\infty)$ and
$w\in A_\infty(\R)$.
The Boyd indices can be computed explicitly for many examples of concrete rearrangement invariant Banach function spaces, see
e.g. \cite[Chapter 4]{BeSh88}.
In particular, we have $p_\E = q_\E = p$ for $\E:=L^{p,q}$ ($1< p<\infty$, $1\leq q\leq \infty$); see \cite[Theorem 4.6]{BeSh88}.

\begin{lemma}\label{wellextr}
Let $\E$ be a rearrangement invariant Banach function space on $(\R, dx)$ such that $1< p_{\E}, q_\E<\infty$. Then the following statements hold.
\begin{itemize}
\item [(i)] For every Muckenhoupt weight $w\in A_{p_\E}(\R)$ there exists a constant $C_{w, \E}$  such that for any family $\cI$ of disjoint bounded intervals in $\R$
\begin{equation}\label{wellextr1}
{C^{-1}_{\E, w}} \|S^\cI f\|_{\E_w}\leq \|S^{\cW^{\cI}}f\|_{\E_w}\leq
{C_{\E, w}} \|S^\cI f\|_{\E_w} 
\end{equation}
and
\begin{equation}\label{FSfs}
  \|M f\|_{\E_w}\leq {C_{\E, w}} \| M^\sharp f\|_{\E_w}
\end{equation}
for every $f\in \E_w$.

\noindent Moreover, if $\cV\subset A_{p_\E}(\R)$ with $\sup_{w\in \cV} [w]_{A_{p_{\E}}}<\infty$, then 
$\sup_{w\in \cV} {C_{\E, w}}<\infty$.

\item [(ii)] For every $r\in (1,\infty)$ and every Muckenhoupt weight
$w\in A_{p_\E}(\R)$ there exists a constant $C_{r,\E,w}$ such that for any family $\cI$ of disjoint intervals in $\R$ 
\begin{equation}\label{wellextr2}
\left\| \left( \sum_{I\in \cI} |S_I f_I|^r \right)^{1/r} \right\|_{\E_w} \leq
C_{r,\E,w} \left\| \left( \sum_{I\in \cI} |f_I|^r \right)^{1/r} \right\|_{\E_w}
\end{equation}
for every $(f_I)_{I\in \cI}\subset \E_w(l^r(\cI))$. 
\end{itemize}
\end{lemma}

The proof follows the idea of the proof of \cite[Lemma 6.3]{RdeF85}, i.e., it is based on the iteration algorithm of the Rubio de Francia extrapolation theory.
We refer the reader to \cite{CrMaPe11} for a recent account of this theory;
 in particular, see the proofs of \cite[Theorems 3.9 and 4.10]{CrMaPe11}.
We provide below main supplementary observations which should be made.

\begin{proof}[{\bf{Proof of Lemma \ref{wellextr}}}]
Note that we can restrict ourself to finite families $\cI$ of disjoint bounded intervals in $\R$. The final estimates obtained below are independent of $\cI$, and a standard limiting argument proves the result in the general case.

According to \cite[Theorem 3.1]{Kurtz80}, for every Muckenhoupt weight $w\in A_2(\R)$ there exists a constant $C_{2,w}$ such that
\begin{equation}\label{2weightest}
C^{-1}_{2,w} \|S^\cI f\|_{L^2(\R,w)} \leq \|S^{\cW^{\cI}}f \|_{L^2(\R,w)} \leq
C_{2,w} \|S^\cI f\|_{L^2(\R,w)}\quad\quad (f\in {L^2(\R,w)}).
\end{equation} Moreover, one can show that $\sup_{w\in\cV}C_{2,w}<\infty$ for every subset $\cV\subset A_2(\R)$ with $\sup_{w\in \cV}[w]_{A_2}<\infty$.

Therefore, we are in a position to adapt the extrapolation techniques from $A_2$ weights; see for example the proof of \cite[Theorem 4.10, p. 76]{CrMaPe11}. Fix $\E$ and $w\in A_{p_\E}(\R)$ as in the assumption. 
Let $\E_{w}'$ be the associate space of $\E_w$, see \cite[Definition 2.3, p. 9]{BeSh88}.
Let  $\cR= \cR_w :\E_{w} \rightarrow \E_{w}$ and $\cR'= \cR'_w :\E_{w}' \rightarrow \E_{w}'$ be defined by
\begin{align*}
\cR h (t) & = \sum_{j=0}^{\infty} \frac{M^j h(t)}{2^j \|M\|^j_{\E_w}}, 
\quad 0\leq h\in \E_{w}, \text{ and} \\
\cR' h (t) & = \sum_{j=0}^{\infty} \frac{S^j h(t)}{2^j \|S\|^j_{\E'_w}}, 
\quad 0\leq h\in \E'_{w},
\end{align*}
where $Sh:= M(h w)/w$ for $h\in \E'_{w}$. As in the proof of 
\cite[Theorem 4.10, p. 76]{CrMaPe11} the following statements are 
easily verified:
\begin{itemize}
\item [(a)]  For every positive $h\in\E_w$ one has 
\begin{align*}
& h\leq \cR h \text{ and } \| \cR h\|_{\E_{w}}\leq 2 \|h\|_{\E_{w}} , 
\text{ and } \\
& \cR h\in A_1 \text{ with }  [\cR h]_{A_1}\leq 2\|M\|_{\E_w} . 
\end{align*}
\item [(b)]  For every positive $h\in\E_w'$ one has 
\begin{align*}
& h\leq \R' h \text{ and } \| \cR' h\|_{\E_{w}'}\leq 2 \|h\|_{\E_{w}'} ,  \text{ and }\\  
& (\cR' h) w\in A_1 \text{ with } [(\cR' h) w]_{A_1}\leq 2\|S\|_{\E_{w}'} .
\end{align*}
\end{itemize} 
The last lines in $(a)$ and $(b)$ follow from the estimates $M(\cR h) \leq 2 \|M\|_{\E_w} \cR h$ and $M((\cR' h)w) \leq 2 \|S\|_{\E'_w} (\cR' h)w)$, respectively, which in turn follow from the definitions of $\cR$ and $\cR'$. 

Note that $f\in L^2(\R, w_{|f|,h})$ for every $f\in \E_w$ and every positive
 $h\in \E'_w$, where
$w_{g,h}:=(\cR g)^{-1}(\cR' h) w$ for every $0\leq g \in \E_w$ and $0\leq h \in \E'_w$. 
 Moreover, by Boyd's interpolation theorem, the Hilbert transform is 
bounded on $\E_w$. Therefore, by the well-known identity relating partial sum operators $S_I$ and the Hilbert transfopozostałyrm, since $\cI$ is finite, we get that $S^\cI f \in \E_w$ for every $f\in \E_w$. 
Similarly, combining Kurtz' inequalities, \cite[Theorem 3.1]{Kurtz80}, with Boyd's interpolation 
theorem, we conclude that $S^{\cW_I} f \in \E_w$ $(I\in \cI)$, and consequently 
$S^{\cW^{\cI}} f \in \E_w$ for every $f\in \E_w$. 

Finally, a close analysis of the proof of \cite[Theorem 4.10]{CrMaPe11}
shows that we can take
$$C_{\E, w}:= 4 \sup\{ C_{2,w_{g,h}}: \, 0\leq g\in \E_w, 0\leq h\in \E'_w, \|g\|_{\E_w}\leq 2, \|h\|_{\E_w'} =1\}.$$

Recall that for every $p\in (1,\infty)$ there exists a constant $C_p>0$ such that  
$\|M\|_{L^p_w} \leq C_{p} [w]^{p'/p}_{A_p}$ for every Muckenhoupt weight $w \in A_p(\R)$; see \cite{Bu93}. A detailed analysis of Boyd's interpolation theorem shows that $\sup_{w\in \cV}\max(\|M\|_{\E_w}, \|S\|_{\E_w}) <\infty$ for every 
$\cV\subset A_{p_{\E}}(\R)$ with $\sup_{w\in \cV}[w]_{A_{p_\E}} <\infty$.
By the so-called reverse factorization (or by H\"older's inequality; see e.g. \cite[Proposition 7.2]{Duo01}), and by properties $(a)$ and $(b)$, we obtain that  $w_{g,h} \in A_2(\R)$ and 
\begin{equation*}
[w_{g,h}]_{A_2}\leq [\cR g]_{A_1}[(\cR' h)w]_{A_1}
\leq 4\|M\|_{\E_w} \|S\|_{\E_{w}'}
\end{equation*} for every $0\leq g \in \E_w$ and $0\leq h \in \E'_w$.
Therefore, on account of the remark on the constants $C_{2,w}$ in \eqref{2weightest}, we get the desired boundedness property of constants $C_{\E, w}$. 
This completes the proof of \eqref{wellextr1}.

Note that, by the weighted Fefferman-Stein inequality, see Remark \ref{upperbound}, and the basic inequality $M^\sharp f\leq 2 Mf$ $(f\in L_{\textrm{loc}}^1(\R))$, the analogous reasoning as before yields \eqref{FSfs}.

For the proof of the part $(ii)$, for fixed $r\in (1,\infty)$ it is sufficient to apply Rubio de Francia's extrapolation algorithm from $A_r$ weights in the same manner as above.
\end{proof}

Let $\cW$ be a well-distributed family of disjoint intervals in $\R$, i.e., there exists $\lambda >1$ such that  
$\sup_{x\in \R} \sum_{I\in \cI} \chi_{\lambda I} (x) <\infty$.
Following \cite[Section 3]{RdeF85}, consider the smooth version of
$S^{\cW}$,  $G=G^{\cW}$, defined as follows:
let $\phi$ be an even, smooth function such that
$\hat \phi (\xi) = 1$
on $\xi\in [-\frac{1}{2}, \frac{1}{2}]$ and $\supp \hat\phi \subset [-\lambda/2,\lambda/2]$.
Let $\phi_I(x):= e^{2\pi i c_I\cdot x} |I|\phi(|I|x)$ $(x\in \R)$,
where $c_I$ stands for the center of an interval $I\in \cW$ and $|I|$ for its length. Then,
$$Gf:=G^{\cW} f : = \left( \sum_{I\in \cW} |\phi_I\star f |^2    \right)^{1/2}\quad (f\in L^2(\R)).$$

Since $\widehat{\phi_I}(\xi)=1$ for  $\xi\in I$, and $\widehat{\phi_I}(\xi)=0$ for
$\xi\notin \lambda I$, by Plancherel's theorem, $G$ is bounded on $L^2(\R)$.

Recall that the crucial step of the proof of \cite[Theorem 6.1]{RdeF85}
consists in showing that the Hilbert space-valued kernel related with $G$ satisfies weak-$(D'_2)$ condition (see \cite[Part IV(E)]{RdeF85} for the definition). This leads to the following pointwise estimates for $G$:
\begin{equation}\label{pointest}
M^\sharp (Gf)(x) \leq C M(|f|^2)(x)^{1/2} \quad \quad (\textrm{ a.e.
} x\in \R)
\end{equation} for every $f\in L^\infty(\R)$ with compact support, and a constant $C$ depending only on $\lambda$. In particular, 
 $G$ is bounded on $L^p(\R, w)$ for every $p>2$ and every Muckenhoupt weight $w\in A_{p/2}(\R)$.

\begin{proof}[{\bf{Proof of Theorem B$(ii)$}}] 
We can assume that $\cI$ is a finite family of bounded intervals in $\R$. By a standard limiting arguments we easily get the general case. 

We start with the proof of the statement of Theorem B$(ii)$ for $q=2$. 
Recall that $p_\E = q_\E = 2$ for $\E:=L^{2,\infty}$; see \cite[Theorem 4.6]{BeSh88}.
Fix $w \in A_1(\R)$ and  $f\in L^\infty(\R)$ 
with compact support. Note that the classical Littlewood-Paley theory shows that  $G^{\cW_I}$ 
is bounded on $L^{2}_w$ for every $I\in \cI$. 
Consequently, $G=G^{\cW^\cI}$ maps $L^2_w$ into itself. 

Therefore, combining Lemma \ref{wellextr}, Lebesgue's differentiation theorem and \eqref{pointest} we get 
\begin{eqnarray*}
 \|S^\cI f\|_{L^{2,\infty}_w} & \leq & C_w \|S^{\cW^\cI} f\|_{L^{2,\infty}_w} \leq C_w \|G f\|_{L^{2,\infty}_w} \leq C_w \|M(G f)\|_{L^{2,\infty}_w}\\
&\leq & C_w \|M^\sharp (G f)\|_{L^{2,\infty}_w} \leq C_w \|M (|f|^2)^{1/2}\|_{L^{2,\infty}_w} 
= C_w \|M (|f|^2)\|_{L^{1,\infty}_w}^{1/2}\\
&\leq& C_w \|f\|_{L^2_w},
\end{eqnarray*}
where $C_w$ is an absolute constant independent on $\cI$ and $f$. The last inequality follows from the fact that the Hardy-Littlewood maximal operator $M$ is of weak $(1,1)$ type. 
Furthermore, one can show that for every subset $\cV \subset A_1(\R)$ with $\sup_{w\in\cV}  [w]_{A_1}<\infty$ we have $\sup_{w\in \cV} C_w<\infty$. 
Since $S^\cI$ is continuous on $L^2_w$ and the space of all functions in $L^\infty(\R)$ with compact support is dense in $L^2_w$  we get the desired boundedness for $S^\cI$.
This completes the proof of the statement of Theorem B$(ii)$ for $q=2$.

We now proceed by interpolation to show the case of $q \in (1,2)$. Let $\cW$ be a well-distributed family of disjoint intervals in $\R$, and $G$ denote the corresponding smooth varsion of $S^\cW$. 
First, it is easily seen that $|\phi_I\star f|\leq (\int \phi dx)Mf$ for every $f\in L_{loc}^1(\R)$ and $I\in \cW$. Moreover, analysis similar to the above shows that 
$G$ maps $L^2_w(\R)$ into $L^{2, \infty}_w(\R)$ for every $w\in A_1(\R)$.
Therefore, for every $w\in A_1(\R)$ the operators 
$$
L^1_w \ni f \mapsto (\phi_I \star f)_{I\in \cW} \in L^{1,\infty}_w(l^\infty)
$$
$$
L_w^2 \ni f \mapsto (\phi_I \star f)_{I\in \cW} \in L^{2,\infty}_w(l^2)
$$ 
are bounded. 
Fix $q\in (1,2)$ and $w\in A_1(\R)$. By interpolation arguments, we conclude that the operator
$$L_w^q \ni f \mapsto (\phi_I \star f)_{I\in \cW} \in L^{q,\infty}_w(l^{q'})$$
is well-defined and bounded.
To show it one can proceed analogously to the proof of a relevant result \cite[Lemma 3.1]{Quek03}. Therefore, we omit details here. 

Since $p_\E = q_\E = q$ for $\E:=L^{q,\infty}$, see \cite[Theorem 4.6]{BeSh88}, by Lemma \ref{wellextr}$(ii)$, for every $w\in A_1(\R)$ we get 
\begin{eqnarray*}
 \|S_q^\cW f\|_{L_w^{q,\infty}} 
& = &  \left\|(\sum_{I\in \cW} |S_I (\phi_I\star f)|^{q'})^{1/q'} \right\|_{L_w^{q,\infty}} \\
&\leq & C_{q,w} \left\|(\sum_{I\in \cW} | \phi_I\star f|^{q'})^{1/q'}\right\|_{L_w^{q,\infty}} \leq C_{q,w} \|f\|_{L_w^q},
\end{eqnarray*}
where $C_{q,w}$ is an absolute constant. 
This completes the proof.
\end{proof}

\begin{remark}
We conclude with the relevant result on $A_2$-weighted $L^2$-estimates for square functions 
$S^\cI$ corresponding to arbitrary families $\cI$ of disjoint intervals in $\R$, i.e., 
$\|S^\cI\|_{2,w} \leq C \|f\|_{2,w}$ $(f\in L_w^2)$.  
According to \cite[Part IV(E)(ii)]{RuRuTo86}, 
these weighted
endpoint estimates can be reached by interpolation provided that $\cI$ is a family such that
$S^\cI$ admits an extension to a bounded operator on (unweighted)
$L^p(\R)$ for some $p<2$. This observation leads to a natural
question: for which partitions $\cI$ of $\R$
 do there exist  {\it local} variants of the Littlewood-Paley decomposition
 theorem, i.e.,  there exists $r\geq 2$ such that  $S^\cI$ is bounded on
$L^p(\R)$ for all $|\frac{1}{p}- \frac{1}{2}|<\frac{1}{r}$.

Recall that L. Carleson, who first noted the possible extension of the classical
Littlewood-Paley inequality for other types of partitions of $\R$, proved
in the special case $\cI:=\{ [n, n+1): n\in \mathbb{Z} \}$ that the corresponding square function $S^\cI$ is bounded on $L^p(\R)$ only if
 $p\geq 2$; see \cite{Carleson67}.
Moreover, it should be noted that such lack of the boundedness of the
square function $S^\cI$ on $L^p(\R)$ for some $p<2$ occurs in
the case of decompositions of $\R$ determined by sequences which
are in a sense not too different from lacunary ones.
Indeed, applying the ideas from \cite[Section 8.5]{EdGa77},
we show below that even in the case of the decomposition $\cI$ of $\R$
determined by a sequence $(a_j)^\infty_{j=0}\subset (0, \infty)$ such that
$a_{j+1} - a_{j}\sim  \lambda^{\phi(j){j}}$, where $\lambda>1$ and
$\phi(j)\rightarrow 0^+$ arbitrary slowly as $j \rightarrow
\infty$, the square function $S^\cI$ is not bounded on $L^p(\R)$ for every $p<2$.

If $I$ is a bounded interval in $\R$, set $f_I$ for the function with 
$\widehat{f_I} = \chi_I$. Then, $|f_I| = \left|\frac{\sin(|I| \pi \cdot)}
{\pi (\cdot)}\right|$, and for every $p>2$ and every $\epsilon >0$ there exists $c>0$ such that
$$\frac{1}{c} |I|^{1/p'} \leq \|f_I\|_p \leq c |I|^{1/p'}$$
for all intervals $I$ with $|I|>\epsilon$. This simply observation
allows to express \cite[Theorem 8.5.4]{EdGa77} for decompositions of $\R$ instead of $\mathbb{Z}$.
Namely, if $a=(a_j)_{j=0}^\infty\subset (0, \infty)$ is an increasing sequence such
 that $a_j - a_{j-1} \rightarrow \infty$ as $j\rightarrow \infty$, and
$\cI_a:=\{ ( - a_0, a_0)\} \cup \{ \pm [a_{j-1}, a_{j})\}_{j\geq 1}$,
then the boundedness of $S^{\cI_a}$ on $L^{p'}(\R)$ for some $p>2$,
 $1/p' + 1/p =1$, implies that there exists a constant $C_p>0$ such that
\begin{equation}\label{equ2}
a_k^{2/p'} \leq C_p \sum_{j=1}^{k}(a_j - a_{j-1})^{2/p'}
\quad (k\geq 1).
\end{equation}
Moreover, it is straightforward to adapt the idea of the proof of
 \cite[Corollary 8.5.5]{EdGa77} to give the following generalization.

Let $a=(a_j)_{j=0}^\infty\subset (0,\infty)$ be an increasing sequence
 such that  $a_{j+1} - a_{j}\sim  \lambda^{\psi(j)}$, where $\lambda>1$, the function
$\psi\in \cC^1([0,\infty))$ is increasing and satisfies the condition:
 $\psi(s)/s \rightarrow 0$ and $\psi '(s) \rightarrow 0$
as $s\rightarrow \infty$.
If the square function $S^{\cI_a}$ were bounded on
 $L^{p'}(\R)$ for some $p>2$, then \eqref{equ2} yields
$$ \left( \int_0^{k-1} \lambda^{\psi(s)} ds\right)^{2/p'}\leq C_p
\int_0^{k+1} \lambda^{\psi(s)2/p'}ds\quad (k\geq 1).
$$
However, this leads to a contradiction with the assumptions on $\psi$.
\end{remark}

\section{Higher dimensional analogue of Theorem A}
The higher dimensional extension of the results due to Coifman, Rubio de Francia and Semmes \cite{CRdeFS88} was established essentially by Q. Xu in \cite{Xu96}; see also M. Lacey \cite[Chapter 4]{Lac07}.

We start with higher dimensional counterparts of some notions from previous
sections. Here and subsequently, we consider only bounded intervals with sides parallel to the axes.

Let $q\geq 1$ and $d\in \N$.
For $h>0$ and  $1\leq k \leq d$ we write $\Delta^{(k)}_h$ for the difference operator, i.e.,
$$\left( \Delta^{(k)}_h m \right) (x):=m(x+ he_k) - m(x)\quad\quad (x\in \R^d)$$
for any function $m:\R^d\rightarrow \C$, where $e_k$ is the $k$-th coordinate vector.
Suppose that $J$ is an interval in $\R^d$ and set
$\overline{J} =:\Pi_{i=1}^{d}[a_i,a_i+h_i]$ with $h_i>0$ $(1\leq i \leq d)$.
We write
$$(\Delta_J m):= \left(\Delta^{(1)}_{h_1}...\Delta^{(d)}_{h_d}m \right) (a),$$
where $a:=(a_1, ..., a_d)$ and $m:\R\rightarrow \C$.
Moreover, for an interval $I$ in $\R^d$ and a function
$m: \R^d \rightarrow \C$ we set
$$\|m\|_{\var_q(I)}:= \sup_{\cJ} \left( \sum_{J\in \cJ} |\Delta_{{J}} m|^q\right)^{1/q},$$
where $\cJ$ ranges over all decompositions of $I$ into subintervals.

Following Q. Xu \cite{Xu96}, see also \cite[Section 4.2]{Lac07},
the spaces $V_q(I)$ for intervals in $\R^d$ are defined inductively as follows.

The definition of $V_q(I)$ $(q\in [1,\infty))$ for one-dimensional intervals is
introduced in Section 1.
Suppose now that $d\in \N\setminus\{1\}$ and fix an interval $I=I_1\times ... \times I_d$ in $\R^d$.
For a function $m:\R^d \rightarrow \C$, we write $m\in V_q(I)$ if
$$\|m\|_{{V_q}(I)}:= \sup_{x\in I}|m(x)| \, +\,
\sup_{x_1\in I_1}
 \|m(x_1, \cdot)\|_{V_q(I_2\times ... \times I_d)}  \, + \, \|m\|_{\var_q(I)} < \infty.$$

Subsequently, $\cD^d$ stands for the family of the dyadic intervals in $\R^d$.
The definition of the spaces $(V_q(\cD^d), \|\cdot\|_{V_q(\cD^d)} )$
 $(d \geq 2,\,\, q\in[1,\infty))$ 
is quite analogous to the corresponding ones in the case of $d=1$ from Section 1.

For a Banach space $X$, an interval $I$ in $\R$ and $q\geq 1$, we consider below the vector-valued variants $V_q(I; X)$, $\cR_q(I; X)$, and  $R_q(I; X)$ of the spaces 
$V(I)$, $\cR_q(I)$, and $R(I)$, respectively. 
  Note that $V_q(I;X)\subset R_p(I;X)$ 
for any $1\leq q <p$ and any interval $I$ in $\R$ with the inclusion norm bounded by a constant depending only on $p$ and $q$; see \cite[Lemma 2]{Xu96}.
Moreover, higher dimensional counterparts of these spaces we define inductively as 
follows: let $I:=\Pi_{i=1}^d I_i$ be a closed interval in $\R^d$ ($d\geq 2)$. 
Set $\widetilde{R}_q(I):= R_q(I_1; \widetilde{R}_q(I_2\times ...\times I_{d}))$ and $\widetilde{V}_q(I) := 
V_q(I_1; \widetilde{V}_q( I_2\times ...\times I_{d}))$, where 
$\widetilde{R}_q(I_d):=R(I_d)$ and $\widetilde{V}_q(I_d):= V_q(I_d)$. 
Recall also that for any $1\leq q <p$ and any interval $I$ in $\R^d$ $(d\geq 1)$
we have
\begin{equation}\label{inclusion}
 V_q(I)\subset \widetilde{V}_q(I) \subset \widetilde{R}_p(I)
\end{equation}
with the inclusion norm bounded by a constant independent of $I$.

Finally, we denote by $A^*_p(\R^d)$ $(p\in
[1,\infty))$ the class of weights on $\R^d$ which satisfy the
strong Muckenhoupt $A_p$ condition.
Note that, in the case of $d=1$, $A^*_p(\R)$ is the classical Muckenhoupt $A_p(\R)$ class ($p\in [1,\infty)$). We refer the reader, e.g., to \cite{Kurtz80} or \cite[Chapter IV.6]{GCRu85}
 for the background on $A^*_p$-weights.

The following complement to \cite[Theorem $(i)$]{Xu96} is the main result of this section.

\begin{thC}\label{ndim}{\it{
  Let $d\geq 2$ and $q\in (1,2]$. Then, $V_q(\cD^d)\subset M_p(\R^d,w)$ for every
$p\geq q$ and every weight $w\in A^*_{p/q}(\R^d)$.\\
$(ii)$ Let $d\geq 2$ and $q>2$. Then, $V_q(\cD^d)\subset M_p(\R^d, w)$ for every
$2\leq p <(\frac{1}{2} - \frac{1}{q})^{-1}$ and every weight $w\in A^*_{p/2}(\R^d)$ with $s_w> (1- p(\frac{1}{2} - \frac{1}{q}))^{-1}$.
}}
\end{thC}

\begin{lemma}\label{main}
 For every $d\in \N$, $q\in(1,2]$, $p> q$, and  every subset $\cV \subset A^*_{p/q}(\R^d)$ with 
$\sup_{w\in \cV}[w]_{A^*_{p/q}(\R^d)} <\infty$ we have $\widetilde{R}_{q}(\cD^d) \subset M_p(\R^d, w)$ $(w\in \cV)$ and 
$$\sup\left\{ \|T_{m\chi_I} \|_{p,w}:  
m\in \widetilde{R}_{q}(\cD^d),\;  \|m\|_{\widetilde{R}_{q}(\cD^d)}\leq 1,\; w\in \cV, \;
I\in \cD^d \right\}
<\infty.$$
\end{lemma}

Here $\widetilde{R}_q(\cD^d)$ $(q\geq 1)$ stands for the space of all functions $m$ defined on $\R^d$ such that $m\chi_I\in \widetilde{R}_q(I)$ for every $I\in\cD^d$ and $\sup_{I\in \cD^d}\|m\chi_I\|_{\widetilde{R}_{q}(I)} < \infty$. Define $\widetilde{V}_q(\cD^d)$ similarly.

The classes $\widetilde{R}_q(\cD^d)$ and $A_p^*(\R^d)$ are well adapted to iterate one-dimensional arguments from the proof of Theorem A$(i)$. 
Therefore, below we give only main supplementary observations should be made.

\begin{proof}[{\bf{Proof of Lemma \ref{main}}}] We proceed by induction on $d$. The proof of the statement of Lemma \ref{main} for $d=1$ and $p=2$ is provided in the proof of Theorem A$(i)$. 
The general case of $d=1$ and $p>q$ follows from this special one by means of Rubio de Francia's extrapolation theorem; see Lemma \ref{extrapol}.

Assume that the statement holds for $d\geq 1$; we will prove it for $d+1$. 
Let $m\in \widetilde{R}_{q}(\cD^{d+1})$ with $\|m\|_{\widetilde{R}_{q}(\cD^{d+1})}\leq 1$. By approximation, we can assume that 
$m_I\in \cR_{q}(I_1; \widetilde{R}_{q}(I_2\times ... \times I_{d+1}) )$ for 
every $I:=I_1\times ... \times I_{d+1} \in \cD^{d+1}$.
Set $m_I :=\sum_{J\in\cI_{I}} \gamma_{I,J} a_{I, J} \chi_J$, where $\gamma_{I,J} \geq 0$ with $\sum_{J}\gamma_{I,J}^{q} \leq 1$ 
and $a_{I, J}\in \widetilde{R}_{q}(I_2\times ... \times I_{d+1})$ with $\|a_{I,J}\|_{ \widetilde{R}_{q}(I_2\times ... \times I_{d+1})} = 1$ for every $I\in \cD^{d+1}$. Here $\cI_{I}$ stands for a decomposition of $I_1$ corresponding to $m_I$.

Let $q\in(1,2]$, $p\geq q'$ and $\cV_{q,p} \subset A^*_{p/q}(\R^{d+1})$ with $\sup_{w\in \cV}[w]_{A^*_{p/q}(\R^{d+1})}<\infty$.  
By Lebesque's differentiation theorem, for every $w\in A^*_{r}(\R^{d+1})$ $(r>1)$ one can easily show that $w(\cdot, y)\in A_{p/q}(\R)$, $w(x,\cdot)\in A_{p/q}^*(\R^{d})$, and   $[w(\cdot,y)]_{A_{p/q}(\R)}, [w(x, \cdot)]_{A_{p/q}^*(\R^{d})}\leq [w]_{A_{p/q}^*(\R^{d+1})}$
for almost every $y\in R^d$ and $x\in \R$;
see e.g. \cite[Lemma 2.2]{Kurtz80}.

Therefore, by induction assumption, for every $q\in(1,2]$ and $  p\in[q',\infty)\setminus \{2\}$  there exists a constant $C_{q,p}>0$ independent of $m$ and $w\in \cV_{q,p}$ such that for every $w \in \cV_{q,p}$:
\begin{equation}\label{est1}
 \sup\left\{ \|T_{a_{I,J}}\|_{p,w(x,\cdot)} : J\in \cI_{I}, \;I\in \cD^{d+1}\right\}\leq 
C_{q,p} 
\quad \quad \textrm{for a.e. } x \in \R.
\end{equation}

Let $f(x,y):=\phi(x)  \rho(y)$ ($(x,y)\in \R^{d+1}$), where $\phi\in S(\R)$ and $\rho\in S(\R^d)$. Note that the set of functions of this form is dense in $L^{q'}(\R^{d+1}, w)$. 
Indeed, by the strong doubling and open ended properties of $A_p^*$-weights, we get
$(1+|\cdot|)^{-dr} w\in L^1(\R^d)$ $(r>1, w\in A_r^*(\R^d))$; see e.g.
\cite[Chapter IX, Proposition 4.5]{Torch04}. Hence, this claim follows from the standard density arguments.  
Moreover, we have  $T_{m_I}f =\sum_J \gamma_{I,J} S_{J}\phi T_{ a_{I,J}}\rho $. 
In the sequel, we consider the case of $q\in (1,2)$ and $q=2$ separately.
For $q\in (1,2)$,
by Fubini's theorem,  we get
\begin{equation*}
 \|T_{m_I} f\|^{q'}_{q', w} \leq \sum_{J\in \cI_I} \int_\R |S_J \phi|^{ q'} \int_{\R^d} |T_{a_{I, J}} \rho|^{q'} w dy dx    \quad \quad (w\in \cV_{q,q'}, I\in \cD^{d+1}).
\end{equation*}
Therefore, by Theorem B$(i)$ and \eqref{est1}, we conclude that
$$\sup\left\{ \|T_{m_I}\|_{q', w}: w\in \cV_{q,q'},\; m\in \widetilde{R}_{q}(\cD^{d+1}), \;
\|m\|_{\widetilde{R}_{q}(\cD^{d+1})}\leq 1, \;I\in \cD^{d+1}\right\}<\infty.$$
Consequently, by Rubio de Francia's extrapolation algorithm, see \cite[Theorem 3]{RdeF84} or \cite[Chapter 3]{CrMaPe11}, the same conclusion holds for all $p>q$. 

For $q=2$, by Fubini's theorem and Minkowski's inequality, we conclude that 
\begin{equation*}
 \|T_{m_I} f\|^{p}_{p, w} \leq  \int_\R |S^{\cI_{I}} \phi(x)|^{p} \left( 
\sum_{J\in \cI_{I}} \gamma^2_{I,J} \|T_{a_{I, J}} \rho\|^{2}_{p, w(x,\cdot)} \right)^{\frac{p}{2}} dx \quad (w\in \cV_{2,p}, I\in \cD^{d+1}).
\end{equation*} for every $p>2$.
Hence, by Theorem \ref{RdeFtheorem} and \eqref{est1}, we get the statement of Lemma \ref{main} also for $q=2$.
\end{proof}

\begin{proof}[{\bf{Proof of Theorem C}}] 
Note first that for every $\cV \subset A_1^*(\R^{d})$ with $\sup_{w\in \cV}[w]_{A_1^*}<\infty$, by the reverse H\"older inequality, there exists $s >1$ such that 
$w^s \in A_1^*(\R^{d})$ $(w\in \cV)$
 and $\sup_{p\geq 2, w\in \cV}[w^s]_{A_{p/2}^*}<\infty$.
Thus,  by Lemma \eqref{main} and an interpolation argument similar to that in the proof of Theorem A$(i)$, we get
$$\sup\left\{ \|T_{m\chi_I}\|_{2, w}: w\in \cV,\; m\in \widetilde{R}_{2}(\cD^{d}), \;
\|m\|_{\widetilde{R}_{2}(\cD^d)}\leq 1, \;I\in \cD^{d}\right\}<\infty .$$

Therefore, as in the proof of Theorem A$(i)$, one can show that for every $q\in(1,2]$
and every subset $\cV \subset A^*_{2/q}(\R^d)$ with 
$N:=\sup_{w\in \cV}[w]_{A^*_{2/q}} <\infty$, there exists a constant 
$\alpha = \alpha(d,q,N)>1$ such that $\widetilde{R}_{\alpha q}(\cD^d) \subset M_2(\R^d, w)$
$(w\in \cV)$ and 
$$\sup\left\{ \|T_{m\chi_I} \|_{2,w}: m\in \widetilde{R}_{\alpha q}(\cD^d),\;  \|m\|_{\widetilde{R}_{\alpha q}(\cD^d)}\leq 1,\; w\in \cV, \; I\in \cD^d \right\} <\infty.$$

Now, by means of \eqref{inclusion}, Kurtz' weighted variant of Littlewood-Paley's inequalities, \cite[Theorem 1]{Kurtz80}, and Rubio de Francia's extrapolation theorem, \cite[Theorem 3]{RdeF84}, the rest of the proof of $(i)$ runs analogously to the corresponding part of the proof of Theorem A$(i)$.

Consequently, by $(i)$, the proof of the part $(ii)$ follows the lines of the proof of Theorem A$(ii)$.
\end{proof}

\end{document}